# The Dodecahedron as a Voronoi Cell

## and its (minor) importance for the Kepler conjecture

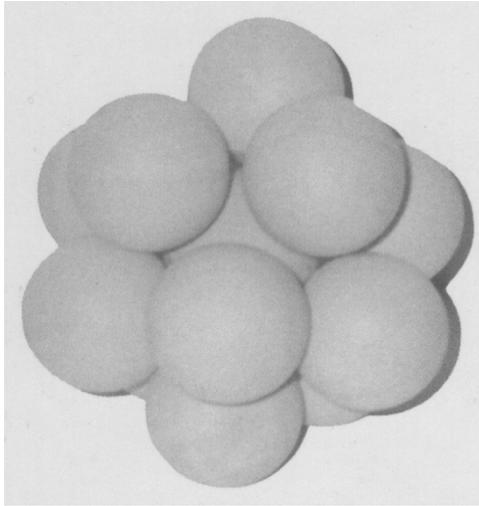

*Icosahedral configuration of 1+12 spheres*
*(This picture is taken from [1] by courtesy of Chuanming Zong.)*

Dr. Max Leppmeier
    Lehrstuhl für Mathematik und ihre Didaktik
    Universität Bayreuth
    max.leppmeier@uni-bayreuth.de

Max Leppmeier studied mathematics and physics at the Ludwig Maximilian University of Munich. He wrote the book "Kugelpackungen von Kepler bis heute" ("Sphere Packings from Kepler until today") and received his doctorate (Ph.D.) from the University of Bayreuth.

**Abstract**
The regular dodecahedron has a 2% smaller volume than the rhombic dodecahedron which is the Voronoi cell of a fcc packing. From this point of view it seems possible that the dodecahedral aspect which is the core of the so-called dodecahedral conjecture, will play a major part for an elementary proof of the Kepler conjecture. In this paper we will show that the icosahedral configuration caused by dodecahedron leads to tetrahedra with significantly larger volume than the fcc fundamental parallelotope tessellation tetrahedra. Therefore on the basis of a tetrahedral based point of view for sphere packing densities we will demonstrate the minor importance of the dodecahedron as a Voronoi cell for the Kepler conjecture.

**Zusammenfassung**
Im Vergleich zum rhombischen Dodekaeder, der Voronoi-Zelle einer fcc-Packung, hat das reguläre Dodekaeder ein um 2% kleineres Volumen. Aus dieser Perspektive erscheint es möglich, dass der dodekaedrische Aspekt, der den Kern der sogenannten dodekaedrischen Vermutung bildet, auch eine entscheidende Rolle für eine elementare Lösung der Kepler-Vermutung spielt. Hier zeigen wir, dass eine mit dem Dodekaeder als Voronoi-Zelle zusammenhängende, ikosaedrische Konfiguration zu Tetraedern mit signifikant größerem Volumen im Vergleich zu den Zerlegungstetraedern eines fcc-Fundamentalparallelotops führt. Somit legt die tetraedrische Perspektive für Kugelpackungsdichten eine untergeordnete Bedeutung des Dodekaeders als Voronoi-Zelle für das Kepler-Hilbert-Problem nahe.

# 1    Introduction

The Kepler conjecture was formulated by Kepler in 1611 and asserts that a densest packing of equal spheres in three-dimensional Euclidean space is given by the „cannonball" packing, or face-centered-cubic (fcc) packing, which fills space with density $\frac{\pi}{\sqrt{18}} \approx 0.74048$ [2, p. 5]. The fcc packing can be described either by its lattice and fundamental parallelotope or by its reciprocal lattice and the Voronoi cell.

The Kepler problem is a famous problem in the field of sphere packings [2] [3] [1] [4] [5] [6]. Hales is regarded as the conqueror of the Kepler problem [7] [8]. The special case of lattice sphere packings in three dimensions is significantly easier and was already solved by Gauß [1] [5].

One peculiarity grounds on the following observation: Around one single sphere we can arrange 12 other spheres not only in a fcc or a related (hcp etc.) configuration (3-6-3), but also in an icosahedral (1-5-5-1) configuration (cf. figure). In the literature we find the dominating conviction that this icosahedral configuration generates a locally larger packing density [9] [1] [2] [10] compared to a fcc configuration. For example, Lagarias formulates:

> "… it is known that an arrangement of 12 unit spheres touching a given unit sphere with their sphere centers being the vertices of a regular dodecahedron[1] yields a Voronoi cell that is a regular dodecahedron of inradius 1, having a ratio of covered to uncovered volume approximately 0.754697, which exceeds $\frac{\pi}{\sqrt{18}} \approx 0.74048$" [2, p. 12].

We find as basic statement the following reason: The Voronoi cell of an icosahedral configuration (which is a regular dodecahedron, cf. the figure in [11, p. 136]) is smaller than the Voronoi cell of a fcc configuration (which is a rhombic dodecahedron, cf. the figure in [11, p. 136]).

Here we focus on tetrahedra and we show that a regular dodecahedron as Voronoi cell leads to icosahedral tetrahedra with a larger volume compared to the volume of a tessellation tetrahedron of a fcc packing. This fact yields a strong argument for the thesis that a regular dodecahedron will not be a relevant Voronoi cell considering the Kepler conjecture.

The topic is related to the famous Newton-Gregory-problem (so-called kissing number problem) [1][5] and to the dodecahedral conjecture [10].

# 2    Definitions and Basics

We base our definitions on [11] [12] [13] [1] [8] [2]. $\mathbb{R}^3$ is the three-dimensional Euclidean space with standard norm $|.|$ . The sphere $K$ is $K = \{x \in \mathbb{R}^3 : |x| \leq 1\}$ with volume $v(K) = \frac{4\pi}{3}$.

Let $X$ be a set of discrete points in $\mathbb{R}^3$ with $|y - x| \geq 2$ for every pair of distinct points $x, y \in X$. Then we call (the Minkowski sum) $X + K = \{x + K; x \in X\}$ a *sphere packing*. If $X$ is a lattice

$$L = \left\{ \sum_1^3 z_i a_i \; ; \; (a_i) \text{ linearly independent}, z_i \in \mathbb{Z} \right\},$$

then we call $L + K$ a *lattice sphere packing*.

In geometrical terms, det $(L)$ is the volume of the *fundamental parallelepiped*

$$FP = \left\{ \sum_1^3 \lambda_i a_i \; ; \; \lambda_i \in [0; 1] \right\}$$

of $L$. A fundamental parallelepiped contains exactly one sphere. Thus the sphere packing density of a given lattice sphere packing can be defined as:

$$\delta(L + K) = \frac{v(K)}{det(L)}$$

Furthermore, the copies $x + FP$ ($x \in L$) yield a space tessellation.

A definition of a non-lattice sphere packing density is more complex [1].

---

[1] Erratum: icosahedron

For a given sphere packing $X + K$ and for $x_i \in X$ we call
$$V(x_i) := \left\{ x;\ |x - x_i| = \min_j |x - x_j| \right\}$$
the *Voronoi cell* of $x_i$. From convex geometry we know that $\bigcup_i V(x_i)$ again yields a tessellation of space. We restrict our consideration on finite Voronoi cells.

## 3   The fcc situation

The common known fcc situation can be described by the following lattice:
$$L_{fcc} = \left\{ z_1 \begin{pmatrix} \sqrt{2} \\ \sqrt{2} \\ 0 \end{pmatrix} + z_2 \begin{pmatrix} \sqrt{2} \\ 0 \\ \sqrt{2} \end{pmatrix} + z_3 \begin{pmatrix} 0 \\ \sqrt{2} \\ \sqrt{2} \end{pmatrix}, z_i \in \mathbb{Z} \right\}$$
In this case the basic vectors span a regular tetrahedron of edge length 2.

The fundamental parallelepiped can be described as
$$FPI := \left\{ \lambda_1 \begin{pmatrix} \sqrt{2} \\ \sqrt{2} \\ 0 \end{pmatrix} + \lambda_2 \begin{pmatrix} \sqrt{2} \\ 0 \\ \sqrt{2} \end{pmatrix} + \lambda_3 \begin{pmatrix} 0 \\ \sqrt{2} \\ \sqrt{2} \end{pmatrix}, \lambda_i \in [0; 1] \right\}$$
expressing the hexagonal character of an fcc packing or as
$$FPII := \left\{ \lambda_1 \begin{pmatrix} 2 \\ 0 \\ 0 \end{pmatrix} + \lambda_2 \begin{pmatrix} 0 \\ 2 \\ 0 \end{pmatrix} + \lambda_3 \begin{pmatrix} 1 \\ 1 \\ \sqrt{2} \end{pmatrix}, \lambda_i \in [0; 1] \right\}$$
formulating the quadratic character of an fcc packing.

Both fundamental parallelepipeds can be tessellated into two regular tetrahedra
$$TIII := \left\{ \lambda \begin{pmatrix} \sqrt{2} \\ \sqrt{2} \\ 0 \end{pmatrix} + \mu \begin{pmatrix} \sqrt{2} \\ 0 \\ \sqrt{2} \end{pmatrix} + \nu \begin{pmatrix} 0 \\ \sqrt{2} \\ \sqrt{2} \end{pmatrix};\ \lambda, \mu, \nu, \lambda + \mu + \nu \in [0; 1] \right\}$$
and four octahedral tetrahedra
$$TIV := \left\{ \lambda \begin{pmatrix} 2 \\ 0 \\ 0 \end{pmatrix} + \mu \begin{pmatrix} 0 \\ 2 \\ 0 \end{pmatrix} + \nu \begin{pmatrix} 1 \\ 1 \\ \sqrt{2} \end{pmatrix};\ \lambda, \mu, \nu, \lambda + \mu + \nu \in [0; 1] \right\}$$
(which develop by dividing an octahedron into four congruent tetrahedra) [14].

Characterizing the different polyhedra we have the

**Lemma.** Let $L_{fcc}$, $FPI$, $FPII$, $TIII$, $TIV$ as defined above. Then the volumes are:
$$v(FPI) = v(FPII) = 4\sqrt{2}$$
and
$$v(TIII) = v(TIV) = \frac{1}{6} v(FPI)$$
**Proof.** The results can be received by elementary geometry. ∎

Here we see that the tessellation tetrahedra of a fundamental parallelepiped have the same volume. This observation will be essential for the following argumentation. (Moreover it can be shown that the packing density within a $TIII$ tetrahedron is slightly larger than the packing density within a $TIV$ tetrahedron [14]. That means: A $TIII$ tetrahedron contains a little bit more than one sixth of a sphere and a $TIV$ tetrahedron contains correspondingly less than one sixth of a sphere.)

Gauß proofed that among all lattice sphere packings the fcc packing is the densest.

**Theorem (Densest lattice sphere packing, Gauß 1831).**
Let $L + K$ be a lattice sphere packing.
Then
$$\sup_L \frac{1}{\det(L)} = \frac{1}{\det(L_{fcc})}$$
and
$$\{L;\ \det(L)\ minimal\} = \{L_{fcc}\}$$

**Proof.** For a detailed proof see for example [15] [5]. ∎

For the fcc configuration packing density it follows from the Lemma:

$$\delta_{fcc} = \frac{v(K)}{v(FPI)} = \frac{\pi}{3\sqrt{2}}$$

We turn to the Voronoi cell of a fcc packing. From the literature we know the following

**Lemma.** (1) Let $L + K$ be a fcc lattice sphere packing and $x_i \in L$. Then each Voronoi cell $V(x_i)$ is a rhombic dodecahedron.
(2) Let $L + K$ be a fcc lattice sphere packing and $V_{fcc}$ a Voronoi cell. Then we have:
$$v(V_{fcc}) = v(FPI) = 4\sqrt{2} \approx 5{,}65685$$
**Proof.** (1) For a proof see for example [15] [1].
(2) Since the Voronoi cells yield a tessellation of $L + K$ as well as the fundamental parallelotopes yield a tessellation of $L + K$ containing exactly one sphere within a Voronoi cell resp. fundamental parallelotope, the volumes of the units must equal. ∎

## 4  Why the dodecahedron could be "better"? – The dodecahedral and the icosahedral situation

We consider the following finite sphere packing configuration: twelve spheres surround and touch a central sphere in an *icosahedral way*. The centers of the surrounding spheres yield the vertices of a regular icosahedron, its center coincides with the center of the central sphere (s. figure).

Comparing the icosahedral configuration with the fcc configuration we see: The surrounding spheres do not touch each other. There is a tiny space between them, and they are a slightly movable. They are ordered in a 1-5-5-1 scheme. In contrast the fcc configuration is fixed; the surrounding spheres do touch each other and are arranged in a 3-6-3 scheme.

We turn to the Voronoi cell of the icosahedral configuration´s central sphere. The facets are determined by the perpendicular bisector planes between the central sphere and a touching sphere. They produce a regular dodecahedron (which is dual to the icosahedron). Its insphere radius corresponds to the radius of the central sphere and thus is 1. Furthermore we have from the formulary literature the

**Lemma.** For a regular dodecahedron with edge length $a$ it holds

for the volume $$v = \frac{a^3}{4}(15 + 7\sqrt{5})$$

and for the insphere radius $$\rho = \frac{a}{2}\sqrt{\frac{25+11\sqrt{5}}{10}}.$$

So we can calculate the volume of a dodecahedral Voronoi cell.

**Lemma.** For the volume of a dodecahedral Voronoi cell it holds:
$$v(V_{dodecahedron}) = \frac{2(15 + 7\sqrt{5})}{\left(\sqrt{\frac{25 + 11\sqrt{5}}{10}}\right)^3} \approx 5{,}55029$$

**Proof.** The insphere radius of the Voronoi cell is $\rho = 1$ since the distance between the central sphere and a touching sphere equals 2.
According to the lemma above we have
$$a = \frac{2\rho}{\sqrt{\frac{25 + 11\sqrt{5}}{10}}} = \frac{2}{\sqrt{\frac{25 + 11\sqrt{5}}{10}}} \approx 0{,}898056$$

and thus for the volume
$$v(V_{dodecahedron}) = \frac{a^3}{4}(15 + 7\sqrt{5})$$

$$v(V_{dodecahedron}) = \frac{\left(\frac{2}{\sqrt{\frac{25 + 11\sqrt{5}}{10}}}\right)^3}{4}(15 + 7\sqrt{5})$$

$$v(V_{dodecahedron}) = \frac{2(15 + 7\sqrt{5})}{\left(\sqrt{\frac{25 + 11\sqrt{5}}{10}}\right)^3} \approx 5{,}55029 \quad \blacksquare$$

This result seems unbelievable: There exists a Voronoi cell with a smaller volume than the volume of the Voronoi cell of a densest lattice sphere packing. L. Fejes Tóth already attached great importance to this peculiarity [9]. Lagarias confirmed the perception as difficulty solving the Kepler conjecture (cf. ch.1).

That perception suggests furthermore a locally higher packing density in the following way:

$$\delta_{dodecahedron} := \frac{v(K)}{v(V_{dodecahedron})} \approx \frac{\frac{4}{3}\pi}{5{,}55029} \approx 75{,}47\%$$

At this point one can have serious doubts whether the fcc packing could be a packing with globally maximum packing density since we have only

$$\delta_{fcc} = \frac{\pi}{3\sqrt{2}} \approx 74{,}05\%.$$

As well, at this point we have at least two counter-arguments: First *space cannot be filled with dodecahedra* as it can be tessellated by rhombic dodecahedra. So one has to calculate with compensation polyhedra which have a locally smaller packing density. (Also a comparison with the planar situation where the inequality of Jessen plays an important role [1] is interesting.) Secondly we should focus on the *initial icosahedral situation*. Here we have tetrahedra as an important basis for packing density considerations.

## 5    The dodecahedron is not ideal! - The solution of the initial icosahedral situation

Within the initial icosahedral situation (cf. fig.) we consider one exemplary tetrahedron defined by one vertex as the center of the central sphere and the other three vertices as centers of three neighbored touching spheres. They yield a regular triangle. We call such a tetrahedron an *icosahedric tetrahedron*.

Calculating the volume of an icosahedric tetrahedron we use the formulary literature.

**Lemma.**  For a regular icosahedron with edge length $a$ it holds:
The circumsphere radius is    $R = \frac{a}{4}\sqrt{10 + 2\sqrt{5}}$.
The insphere radius is    $\rho = \frac{a}{12}\sqrt{3}(3 + \sqrt{5})$. $\blacksquare$

Since the centers of the touching spheres yield the vertices of the icosahedron its circumsphere radius is exactly 2 and we conclude the

**Lemma.**  Within the icosahedral configuration the distance between two neighbored touching spheres is

$$a = \frac{8}{\sqrt{10 + 2\sqrt{5}}} \approx 2{,}1029.$$

**Proof.**  From $R = 2$ we have by the Lemma above

$$\frac{a}{4}\sqrt{10 + 2\sqrt{5}} = 2$$

and thus the assumption. $\blacksquare$

This result quantifies the space between the touching spheres. Within a range of 5% related to the radius they are slightly movable on the central sphere.

Now we calculate the icosahedral tetrahedron. Its height corresponds to the insphere radius of the icosahedron. We have the

**Lemma.**  For the icosahedral tetrahedron it holds:
The height is    $h = \rho = \frac{2\sqrt{3}(3+\sqrt{5})}{3\sqrt{10+2\sqrt{5}}} \approx 1{,}58931$.
**Proof.**  We have for the base edge length of the icosahedral tetrahedron

$$a = \frac{8}{\sqrt{10 + 2\sqrt{5}}}$$

and thus for the insphere radius

$$\rho = \frac{a}{12}\sqrt{3}(3 + \sqrt{5}) = \frac{8}{12\sqrt{10 + 2\sqrt{5}}}\sqrt{3}(3 + \sqrt{5})$$

$$\rho = \frac{2\sqrt{3}(3 + \sqrt{5})}{3\sqrt{10 + 2\sqrt{5}}} \approx 1{,}58931$$

confirming the claim. ∎

Here we see a quite differentiated pattern: The movability and the larger distance of the touching spheres produce a longer edge length and a larger base area of the icosahedral tetrahedron. The essential question is therefore: can the smaller height compensate this effect – compared to a regular tetrahedron?
We find the answer in the following

**Theorem.**
  (1) A regular tetrahedron with edge length $a = 2$ has
      base area $A = \sqrt{3} \approx 1{,}73205$, height $h = \frac{2}{3}\sqrt{6} \approx 1{,}63299$
      and volume $V = \frac{2}{3}\sqrt{2} \approx 0{,}94281$.
  (2) The icosahedral tetrahedron with base edge length $a = \frac{8}{\sqrt{10+2\sqrt{5}}} \approx 2{,}1029$ has
      base area $A = \frac{\sqrt{3}}{4}\left(\frac{8}{\sqrt{10+2\sqrt{5}}}\right)^2 \approx 1{,}91491$, height $h = \frac{2\sqrt{3}(3+\sqrt{5})}{3\sqrt{10+2\sqrt{5}}} \approx 1{,}58931$
      and volume $V = \frac{1}{3}\frac{\sqrt{3}}{4}\left(\frac{8}{\sqrt{10+2\sqrt{5}}}\right)^2 \frac{2\sqrt{3}(3+\sqrt{5})}{3\sqrt{10+2\sqrt{5}}} \approx 1{,}01446$.

**Proof.** We get the proof from the lemmas above and elementary calculations. ∎

The answer is clear: No.
The base area of an icosahedral tetrahedron is more than 10% larger than the base area of a regular tetrahedron: The height is not even 3% smaller than the height of a regular tetrahedron and thus cannot compensate the base area deficits.

## 6   Interpretation

Lattices, fundamental parallelotopes and tessellation tetrahedra are essential parts of the Gauß proof. The volume of a tetrahedron is minimal (under certain conditions) for a regular tetrahedron $TIII$ or for an octahedral tetrahedron $TIV$. Both yield the tessellation tetrahedra of a fcc packing configuration.

   Within a fcc sphere packing we have six tessellation tetrahedra for each fundamental parallelotope ($FPI$ or $FPII$): two regular tetrahedra and four octahedral tetrahedra. Both types have the same volume, one sixth of the fundamental parallelotope´s volume.

   The fcc packing can be seen as packing of tessellation tetrahedra with volume $V = \frac{2}{3}\sqrt{2} \approx 0{,}94281$. In contrast, an icosahedral tetrahedron has the volume $V = \frac{1}{3}\frac{\sqrt{3}}{4}\left(\frac{8}{\sqrt{10+2\sqrt{5}}}\right)^2 \frac{2\sqrt{3}(3+\sqrt{5})}{3\sqrt{10+2\sqrt{5}}} \approx 1{,}01446$.

   Besides a detailed consideration of the packing densities within the tetrahedra (what is given for the regular and the octahedral tetrahedra in [14]) we can assert: The volume of an icosahedral tetrahedron evidences in a strong way that an icosahedral configuration will not play a major or even crucial role in the context of the Kepler conjecture. Thus the regular dodecahedron is a detour and not essential for an elementary proof of the Kepler conjecture.

## 7   Sideview to Solid State Physics

We know from solid state physics the meaning of the reciprocal lattice and the Brillouin zone [3].

   For example, the Brillouin zone of a simple cubic packing (sc) again is a cube, but the Brillouin zone of a body centered cubic packing (bcc) is a rhombic dodecahedron i.e. the Voronoi cell of a fcc packing. However nobody claims to rethink whether the bcc packing would be denser than the fcc packing.

Again from solid state physics we know that the Wigner-Seitz cell (corresponding to the Voronoi cell) contains the same volume as the primitive cell and also contains one sphere as the primitive cell. For these reasons we have two kinds of tessellation with the same packing density within a Wigner-Seitz cell or a primitive cell.

If the tessellation aspect drops out, we can only consider the local situation. As one consequence, the volumes of a local Wigner-Seitz cell and of a local primitive cell differ from each other: One must be larger than the other.

If we had a voluminous Wigner-Seitz cell what causes a less voluminous primitive cell, we would have a serious problem in respect of the minimal volume property of a regular tetrahedron and the Gauß proof.

But we just discussed a small volume of the dodecahedron as a local Wigner-Seitz cell. This situation causes a voluminous primitive cell and is consistent with the minimal volume property of a regular tetrahedron and the Gauß proof, and therefore constitutes no serious problem.

**Acknowledgements**
I owe many thanks to V. Ulm for his support, A. Beutelspacher for his encouragement, J.M. Wills for a long-time scientific mentoring and P.M. Gruber for a great summer school on the theme of sphere packings and convex geometry in La Villa in 2000.